\documentclass[psamsfonts,english]{amsart}

    \usepackage[T1]{fontenc}
    \usepackage{babel}
    \usepackage{natbib}
    \usepackage{amsmath}
    \usepackage{amssymb}
    \usepackage{amsthm}
    \usepackage[mathcal]{eucal}
    \usepackage{url}

    \newtheorem{teo}{Theorem}
    \newtheorem{lem}[teo]{Lemma}
    \newtheorem{kor}[teo]{Corollary}
    \theoremstyle{definition}
    
    \newtheorem{alg}{Algorithm}
    \theoremstyle{remark}
    \newtheorem{rem}[teo]{Remark}

    \newcommand{\FF}{\mathbb{F}}
    \newcommand{\ZZ}{\mathbb{Z}}
    \newcommand{\PP}{\mathbb{P}}
    \newcommand{\jac}[1]{\mathcal{J}_{#1}}
    \newcommand{\grp}[1]{\langle #1 \rangle}
    \newcommand{\gruppen}{\Gamma}

    \DeclareMathOperator{\Div}{Div}
    \DeclareMathOperator{\divisor}{div}

\begin{document}

\title[Generators of Jacobians of Hyperelliptic Curves]
{Generators of Jacobians of Hyperelliptic Curves}

\author[C.R. Ravnshøj]{Christian Robenhagen Ravnshøj}


\address{Department of Mathematical Sciences \\
University of Aarhus \\
Ny Munkegade \\
Building 1530 \\
DK-8000 Aarhus C}

\email{cr@imf.au.dk}

\thanks{Research supported in part by a Ph.D. grant from CRYPTOMAThIC}




\begin{abstract}
This paper provides a probabilistic algorithm to determine generators of the $m$-torsion subgroup
of the Jacobian of a hyper\-elliptic curve of genus two.
\end{abstract}

\maketitle

\section{Introduction}\label{sec:intro}

Let $C$ be a hyperelliptic curve of genus two defined over a prime field $\FF_p$, and~$\jac{C}$ the Jacobian of $C$.
Consider the rational subgroup $\jac{C}(\FF_p)$. $\jac{C}(\FF_p)$ is a finite abelian group, and
    $$\jac{C}(\FF_p)\simeq\ZZ/n_1\ZZ\oplus\ZZ/n_2\ZZ\oplus\ZZ/n_3\ZZ\oplus\ZZ/n_4\ZZ,$$
where $n_i\mid n_{i+1}$ and $n_2\mid p-1$. \cite{frey-ruck} shows that if $m\mid p-1$, then the discrete logarithm
problem in the rational $m$-torsion subgroup $\jac{C}(\FF_p)[m]$ of $\jac{C}(\FF_p)$ can be reduced to the
corresponding problem in $\FF_p^\times$ \cite[corollary~1]{frey-ruck}. In the proof of this result it is claimed that
the non-degeneracy of the Tate pairing can be used to determine whether $r$ random elements of the finite group
$\jac{C}(\FF_p)[m]$ in fact is an independent set of generators of $\jac{C}(\FF_p)[m]$. This paper provides an
explicit, probabilistic algorithm to determine generators of $\jac{C}(\FF_p)[m]$.

In short, the algorithm outputs elements $\gamma_i$ of the Sylow-$\ell$ subgroup $\gruppen_\ell$ of the rational subgroup $\gruppen=\jac{C}(\FF_p)$, such that
$\gruppen_\ell=\bigoplus_i\grp{\gamma_i}$ in the following steps:
\begin{enumerate}
 \item Choose random elements $\gamma_i\in\gruppen_\ell$ and $h_j\in\jac{C}(\FF_p)$, $i,j\in\{1,\dots,4\}$.\label{choose:intro}
 \item Use the non-degeneracy of the tame Tate pairing $\tau$ to \emph{diagonalize} the sets $\{\gamma_i\}_i$ and
 $\{h_j\}_j$ with respect to $\tau$; i.e. modify the sets such that $\tau(\gamma_i,h_j)=1$ if $i\neq j$ and
 $\tau(\gamma_i,h_i)$ is an $\ell^{\mathrm{th}}$ root of unity.\label{step:diagonaliser}
 \item If $\prod_i|\gamma_i|<|\gruppen_\ell|$ then go to step~\ref{choose:intro}.
 \item Output the elements $\gamma_1$, $\gamma_2$, $\gamma_3$ and $\gamma_4$.
\end{enumerate}
The key ingredient of the algorithm is the diagonalization in step~\ref{step:diagonaliser}; this process will be
explained in section~\ref{sec:generators}.

We will write $\grp{\gamma_i|i\in I}=\grp{\gamma_i}_i$ and $\bigoplus_{i\in I}\grp{\gamma_i}=\bigoplus_i\grp{\gamma_i}$ if the index set
$I$ is clear from the context.

\section{Hyperelliptic curves}

A hyperelliptic curve is a smooth, projective curve $C\subseteq\PP^n$ of genus at least two with a separable, degree
two morphism $\phi:C\to\PP^1$. In the rest of this paper, let~$C$ be a hyperelliptic curve of genus two defined over a
prime field $\FF_p$ of characteristic~$p>2$. By the Riemann-Roch theorem there exists an embedding $\psi:C\to\PP^2$,
mapping $C$ to a curve given by an equation of the form
    $$y^2=f(x),$$
where $f\in\FF_p[x]$ is of degree six and have no multiple roots \cite[see][chapter~1]{cassels}.

The set of principal divisors $\mathcal{P}(C)$ on $C$ constitutes a subgroup of the degree zero divisors $\Div_0(C)$. The
Jacobian $\jac{C}$ of $C$ is defined as the quotient
    $$\jac{C}=\Div_0(C)/\mathcal{P}(C).$$
Consider the subgroup $\jac{C}(\FF_p)<\jac{C}$ of $\FF_p$-rational elements. There exist numbers~$n_i$, such that
    \begin{equation}\label{eq:rank4}
    \jac{C}(\FF_p)\simeq\ZZ/n_1\ZZ\oplus\ZZ/n_2\ZZ\oplus\ZZ/n_3\ZZ\oplus\ZZ/n_4\ZZ,
    \end{equation}
where $n_i\mid n_{i+1}$ and $n_2\mid p-1$ \cite[see][proposition~5.78, p.~111]{hhec}. We wish to determine generators
of the $m$-torsion subgroup $\jac{C}(\FF_p)[m]<\jac{C}(\FF_p)$, where $m\mid |\jac{C}(\FF_p)|$ is the largest number
such that $\ell\mid p-1$ for every prime number $\ell\mid m$.

\section{Finite abelian groups}

\cite{miller} shows the following theorem.

\begin{teo}\label{teo:ssh-generator}
Let $G$ be a finite abelian group of torsion rank $r$. Then for $s\geq r$ the probability that a random $s$-tuple of
elements of $G$ generates $G$ is at least
    $$\frac{C_r}{\log\log|G|}$$
if $s=r$, and at least $C_s$ if $s>r$, where $C_s>0$ is a constant depending only on $s$ (and not on $|G|$).
\end{teo}

\begin{proof}
\cite[theorem~3, p.~251]{miller}
\end{proof}

Combining theorem~\ref{teo:ssh-generator} and equation~\eqref{eq:rank4}, we expect to find generators of $\gruppen[m]$
by choosing $4$ random elements $\gamma_i\in\gruppen[m]$ in approximately $\frac{\log\log |\gruppen[m]|}{C_4}$
attempts.

To determine whether the generators are independent, i.e. if $\grp{\gamma_i}_i=\bigoplus_i\grp{\gamma_i}$, we need to
know the subgroups of a cyclic $\ell$-group $G$. These are determined uniquely by the order of $G$, since
    $$\{0\}<\grp{\ell^{n-1}g}<\grp{\ell^{n-2}g}<\dots<\grp{\ell g}<G$$
are the subgroups of the group $G=\grp{g}$ of order $\ell^n$. The following corollary is an immediate consequence of
this observation.

\begin{kor}\label{kor:UGsnit}
Let $U_1$ and $U_2$ be cyclic subgroups of a finite group $G$. Assume $U_1$ and $U_2$ are $\ell$-groups. Let
$\grp{u_i}<U_i$ be the subgroups of order $\ell$. Then
    $$U_1\cap U_2=\{e\}\Longleftrightarrow\grp{u_1}\cap\grp{u_2}=\{e\}.$$
Here $e\in G$ is the neutral element.
\end{kor}

\section{The tame Tate pairing}

Let $\gruppen=\jac{C}(\FF_p)$ be the rational subgroup of the Jacobian. Consider a number
$\lambda\mid\gcd(|\gruppen|,p-1)$. Let $g\in\gruppen[\lambda]$ and $h=\sum_ia_i P_i\in\gruppen$ be divisors with no
points in common, and let
    $$\overline{h}\in\gruppen/\lambda\gruppen$$
denote the class containing the divisor~$h$. Furthermore, let $f\in\FF_{p}(C)$ be a rational function on $C$ with
divisor $\divisor(f)=\lambda g$. Set $f(h)=\prod_if(P_i)^{a_i}$. Then
    $$e_\lambda(g,\overline{h})=f(h)$$
is a well-defined pairing
$\gruppen[\lambda]\times\gruppen/\lambda\gruppen\longrightarrow\FF_{p}^\times/(\FF_{p}^\times)^\lambda$, the \emph{Tate
pairing}; cf.~\cite{galbraith}. Raising to the power $\frac{p-1}{\lambda}$ gives a well-defined element in the subgroup
$\mu_\lambda<\FF_{p}^\times$ of the $\lambda^{\mathrm{th}}$ roots of unity. This pairing
    $$\tau_\lambda:\gruppen[\lambda]\times\gruppen/\lambda\gruppen\longrightarrow\mu_\lambda$$
is called the \emph{tame Tate pairing}.

Since the class $\overline{h}$ is represented by the element $h\in\gruppen$, we will write $\tau_\lambda(g,h)$ instead of
$\tau_\lambda(g,\overline{h})$. Furthermore, we will omit the subscript $\lambda$ and just write $\tau(g,h)$, since the
value of $\lambda$ will be clear from the context.

\cite{hess} gives a short and elementary proof of the following theorem.

\begin{teo}\label{teo:tatepairing}
The tame Tate pairing $\tau$ is bilinear and non-degenerate.
\end{teo}

\begin{kor}\label{kor:tatepairing}
For every element $g\in\gruppen$ of order $\lambda$ an element $h\in\gruppen$ exists, such that
$\mu_\lambda=\grp{\tau(g,h)}$.
\end{kor}

\begin{proof}
\cite[corollary~8.1.1., p.~98]{sil} gives a similar result for elliptic curves and the Weil pairing. The proof of this
result only uses that the pairing is bilinear and non-degenerate. Hence it applies to corollary~\ref{kor:tatepairing}.
\end{proof}

\begin{rem}
In the following we only need the existence of the element $h\in\gruppen$, such that $\mu_\lambda=\grp{\tau(g,h)}$; we
do not need to find it.
\end{rem}

\section{Generators of $\gruppen[m]$}\label{sec:generators}

As in the previous section, let $\gruppen=\jac{C}(\FF_p)$ be the rational subgroup of the Jacobian. We are searching
for elements $\gamma_i\in\gruppen[m]$ such that $\gruppen[m]=\bigoplus_i\grp{\gamma_i}$. As an abelian group,
$\gruppen[m]$ is the direct sum of its Sylow subgroups. Hence, we only need to find generators of the Sylow subgroups of
$\gruppen[m]$.

Set $N=|\gruppen|$ and let $\ell\mid\gcd(N,p-1)$ be a prime number. Choose four random elements~$\gamma_i\in\gruppen$. Let
$\gruppen_\ell<\gruppen$ be the Sylow-$\ell$ subgroup of $\gruppen$, and set $N_\ell=|\gruppen_\ell|$. Then
$\frac{N}{N_\ell}\gamma_i\in \gruppen_\ell$. Hence, we may assume that $\gamma_i\in \gruppen_\ell$. If all the elements $\gamma_i$ are equal to zero,
then we choose other elements $\gamma_i\in\gruppen$. Hence, we may assume that some of the elements~$\gamma_i$ are non-zero.

Let $|\gamma_i|=\lambda_i$, and re-enumerate the $\gamma_i$'s such that $\lambda_i\leq\lambda_{i+1}$. Since some of the $\gamma_i$'s are non-zero, we may choose an index $\nu\leq 4$, such that $\lambda_\nu\neq 1$ and $\lambda_i=1$ for $i<\nu$. Choose~$\lambda_0$ minimal such that $\lambda=\frac{\lambda_\nu}{\lambda_0}\mid p-1$. Then $\FF_p$ contains an
element $\zeta$ of order $\lambda$. Now set $g_i=\frac{\lambda_i}{\lambda}\gamma_i$, $\nu\leq i\leq 4$. Then
$g_i\in\gruppen[\lambda]$, $\nu\leq i\leq 4$. Finally, choose four random elements $h_i\in\gruppen$.

Let
    $$\tau:\gruppen[\lambda]\times\gruppen/\lambda\gruppen\longrightarrow\grp{\zeta}$$
be the tame Tate pairing. Define remainders $\alpha_{ij}$ modulo $\lambda$ by
    $$\tau(g_i,h_j)=\zeta^{\alpha_{ij}}.$$
By corollary~\ref{kor:tatepairing}, for any of the elements $g_i$ we can choose an element $h\in\gruppen$, such that
$|\tau(g_i,h)|=\lambda$. Assume that
$\gruppen/\lambda\gruppen=\grp{\overline{h}_1,\overline{h}_2,\overline{h}_3,\overline{h}_4}$. Then
$\overline{h}=\sum_iq_i\overline{h}_i$, and so
    $$\tau(g_i,h)=\zeta^{\alpha_{i1}q_1+\alpha_{i2}q_2+\alpha_{i3}q_3+\alpha_{i4}q_4}.$$
If $\alpha_{ij}\equiv 0\pmod{\ell}$, $1\leq j\leq 4$, then $|\tau(g_i,h)|<\lambda$. Hence, if
$\gruppen/\lambda\gruppen=\grp{\overline{h}_1,\overline{h}_2,\overline{h}_3,\overline{h}_4}$, then for all
$i\in\{\nu,\dots,4\}$ we can choose a $j\in\{1,\dots,4\}$, such that $\alpha_{ij}\not\equiv 0\pmod{\ell}$.

Enumerate the $h_i$ such that $\alpha_{44}\not\equiv 0\pmod{\ell}$. Now assume a number $j<4$ exists, such that
$\alpha_{4j}\not\equiv 0\pmod{\lambda}$. Then $\zeta^{\alpha_{4j}}=\zeta^{\beta_1\alpha_{44}}$, and replacing $h_j$
with $h_j-\beta_1h_4$ gives $\alpha_{4j}\equiv 0\pmod{\lambda}$. So we may assume that
    $$\alpha_{41}\equiv\alpha_{42}\equiv\alpha_{43}\equiv 0\pmod{\lambda}\qquad\textrm{and}\qquad\alpha_{44}\not\equiv 0\pmod{\ell}.$$
Assume similarly that a number $j<4$ exists, such that $\alpha_{j4}\not\equiv 0\pmod{\lambda}$. Now
set~$\beta_2\equiv\alpha_{44}^{-1}\alpha_{j4}\pmod{\lambda}$. Then $\tau(g_j-\beta_2g_4,h_4)=1$. So we may also assume
that
    $$\alpha_{14}\equiv\alpha_{24}\equiv\alpha_{34}\equiv 0\pmod{\lambda}.$$
Repeating this process recursively, we may assume that
    $$\alpha_{ij}\equiv 0\pmod{\lambda}\qquad\textrm{and}\qquad\alpha_{44}\not\equiv 0\pmod{\ell}.$$
Again $\nu\leq i\leq 4$ and $1\leq j\leq 4$.

The discussion above is formalized in the following algorithm.

\begin{alg}\label{alg:1}
As input we are given a hyperelliptic curve $C$ of genus two defined over a prime field $\FF_p$, the number
$N=|\gruppen|$ of $\FF_p$-rational elements of the Jacobian, and a prime factor $\ell\mid\gcd(N,p-1)$. The algorithm
outputs elements $\gamma_i\in \gruppen_\ell$ of the Sylow-$\ell$ subgroup $\gruppen_\ell$ of $\gruppen$, such that
$\grp{\gamma_i}_i=\bigoplus_i\grp{\gamma_i}$ in the following steps.
\begin{enumerate}
 \item Compute the order $N_\ell$ of the Sylow-$\ell$ subgroup of $\gruppen$.
 \item Choose elements $\gamma_i\in\gruppen$, $i\in I:=\{1,2,3,4\}$. Set $\gamma_i:=\frac{N}{N_\ell}\gamma_i$. \label{step:choose-x}
 \item Choose elements $h_j\in\gruppen$, $j\in J:=\{1,2,3,4\}$. \label{step:choose-h}
 \item Set $K:=\{1,2,3,4\}$.
 \item For $k'$ from $0$ to $3$ do the following:
 \begin{enumerate}
  \item Set $k:=4-k'$.
  \item If $\gamma_i=0$, then set $I:=I\setminus\{i\}$. If $|I|=0$, then go to step~\ref{step:choose-x}.
  \item Compute the orders $\lambda_\kappa:=|\gamma_\kappa|$, $\kappa\in K$.
        Re-enumerate the $\gamma_\kappa$'s such that $\lambda_\kappa\leq \lambda_{\kappa+1}$, $\kappa\in K$.
        Set $I:=\{5-|I|,6-|I|,\dots,4\}$.
  \item Set $\nu:=\min(I)$, and choose $\lambda_0$ minimal such that
        $\lambda:=\frac{\lambda_\nu}{\lambda_0}\mid p-1$.
        Set $g_\kappa:=\frac{\lambda_\kappa}{\lambda}\gamma_\kappa$, $\kappa\in I\cap K$.
  \begin{enumerate}
   \item If $g_k=0$, then go to step~\ref{step:last}.
   \item If $\tau(g_k,h_j)^{\lambda/\ell}=1$ for all $j\leq k$, then go to step~\ref{step:choose-h}.
  \end{enumerate}
  \item Choose a primitive $\lambda^{\mathrm{th}}$ root of unity $\zeta\in\FF_p$.
        Compute $\alpha_{kj}$ and $\alpha_{\kappa k}$ from
        $\tau(g_k,h_j)=\zeta^{\alpha_{kj}}$ and $\tau(g_\kappa,h_k)=\zeta^{\alpha_{\kappa k}}$, $1\leq j<k$, $\kappa\in I\cap K$.
        Re-enumerate $h_1,\dots,h_k$ such that $\alpha_{kk}\not\equiv 0\pmod{\ell}$.
  \item For $1\leq j<k$, set $\beta\equiv\alpha_{kk}^{-1}\alpha_{kj}\pmod{\lambda}$ and $h_j:=h_j-\beta h_k$.
  \item For $\kappa\in I\cap K\setminus\{k\}$, set $\beta\equiv\alpha_{kk}^{-1}\alpha_{\kappa k}\pmod{\lambda}$
        and $\gamma_\kappa:=\gamma_\kappa-\beta\frac{\lambda_k}{\lambda_\kappa}\gamma_k$.
  \item Set $K:=K\setminus\{k\}$.
 \end{enumerate}
 \item Output $\gamma_1$, $\gamma_2$, $\gamma_3$ and $\gamma_4$.\label{step:last}
\end{enumerate}
\end{alg}

{\samepage
\begin{rem}\label{rem:runningtime}
Algorithm~\ref{alg:1} consists of a small number of
\begin{enumerate}
 \item calculations of orders of elements $\gamma\in\gruppen_\ell$,\label{item:orden}
 \item multiplications of elements $\gamma\in\gruppen$ with numbers $a\in\ZZ$,\label{item:gange}
 \item additions of elements $\gamma_1,\gamma_2\in\gruppen$,\label{item:addition}
 \item evaluations of pairings of elements $\gamma_1,\gamma_2\in\gruppen$ and\label{item:parring}
 \item solving the discrete logarithm problem in $\FF_p$, i.e. to determine $\alpha$ from $\zeta$ and $\xi=\zeta^\alpha$.\label{item:DL}
\end{enumerate}
By \cite[proposition~9]{miller}, the order $|\gamma|$ of an element $\gamma\in\gruppen_\ell$ can be calculated in time
$O(\log^3 N_\ell)\mathcal{A}_\gruppen$, where $\mathcal{A}_\gruppen$ is the time for adding two elements of $\gruppen$.
A multiple $a\gamma$ or a sum $\gamma_1+\gamma_2$ is computed in time $O(\mathcal{A}_\gruppen)$. By \cite{frey-ruck},
the pairing $\tau(\gamma_1,\gamma_2)$ of two elements $\gamma_1,\gamma_2\in\gruppen$ can be evaluated in time $O(\log N_\ell)$. Finally, by
\cite{pohlig-hellmann} the discrete logarithm problem in $\FF_p$ can be solved in time $O(\log p)$. We may assume that
addition in $\gruppen$ is easy, i.e. that $\mathcal{A}_\gruppen<O(\log p)$. Hence algorithm~\ref{alg:1} runs in
expected time $O(\log p)$.
\end{rem}}

{\samepage Careful examination of algorithm~\ref{alg:1} gives the following lemma.

\begin{lem}\label{lem:diagonal}
Let $\gruppen_\ell$ be the Sylow-$\ell$ subgroup of $\gruppen$, $\ell\mid p-1$. Algorithm~\ref{alg:1} determines
elements $\gamma_i\in \gruppen_\ell$ and $h_i\in\gruppen$, $1\leq i\leq 4$, such that one of the following cases holds.
\begin{enumerate}
 \item $\alpha_{11}\alpha_{22}\alpha_{33}\alpha_{44}\not\equiv 0\pmod{\ell}$ and $\alpha_{ij}\equiv 0\pmod{\lambda}$, $i\neq
 j$, $i,j\in\{1,2,3,4\}$.\label{case:dia1}
 \item $\gamma_1=0$, $\alpha_{22}\alpha_{33}\alpha_{44}\not\equiv 0\pmod{\ell}$ and $\alpha_{ij}\equiv 0\pmod{\lambda}$, $i\neq
 j$, $i,j\in\{2,3,4\}$.
 \item $\gamma_1=\gamma_2=0$, $\alpha_{33}\alpha_{44}\not\equiv 0\pmod{\ell}$ and $\alpha_{ij}\equiv 0\pmod{\lambda}$, $i\neq
 j$, $i,j\in\{3,4\}$. 
 \item $\gamma_1=\gamma_2=\gamma_3=0$.
\end{enumerate}
If $|\gamma_i|=\lambda_i$, then $\lambda_i\leq \lambda_{i+1}$. Set $\nu=\min\{i|\lambda_i\neq 1\}$, and define
$\lambda_0$ as the least number, such that $\lambda=\frac{\lambda_\nu}{\lambda_0}\mid p-1$. Set
$g_i=\frac{\lambda_i}{\lambda}\gamma_i$, $\nu\leq i\leq 4$. Then the numbers $\alpha_{ij}$ above are determined by
    $$\tau(g_i,h_j)=\zeta^{\alpha_{ij}},$$
where $\tau$ is the tame Tate pairing $\gruppen[\lambda]\times\gruppen/\lambda\gruppen\to\mu_\lambda=\grp{\zeta}$.
\end{lem}
}

\begin{teo}\label{teo:p-1}
Algorithm~\ref{alg:1} determines elements $\gamma_1$, $\gamma_2$, $\gamma_3$ and $\gamma_4$ of the Sylow-$\ell$
subgroup of $\gruppen$, $\ell\mid p-1$, such that $\grp{\gamma_i}_i=\bigoplus_i\grp{\gamma_i}$.
\end{teo}

\begin{proof}
Choose elements $\gamma_i,h_i\in\gruppen$ such that the conditions of lemma~\ref{lem:diagonal} are fulfilled. Set $\lambda_i=|\gamma_i|$,
and let $\nu=\min\{i|\lambda_i\neq 1\}$. Define $\lambda_0$ as the least number, such that
$\lambda=\frac{\lambda_\nu}{\lambda_0}\mid p-1$. Set $g_i=\frac{\lambda_i}{\lambda}\gamma_i$. Then the $\alpha_{ij}$'s
from lemma~\ref{lem:diagonal} are determined by
    $$\tau(g_i,h_j)=\zeta^{\alpha_{ij}}.$$
We only consider case~\ref{case:dia1} of lemma~\ref{lem:diagonal}, since the other cases follow similarly. We start by
determining $\grp{\gamma_3}\cap\grp{\gamma_4}$. Assume that $g_3=ag_4$. Then
    $$1=\tau(g_3,h_4)=\tau(ag_4,h_4)=\zeta^{a\alpha_{44}},$$
i.e. $a\equiv 0\pmod{\lambda}$. Hence $\grp{\gamma_3}\cap\grp{\gamma_4}=\{0\}$. Then we determine
$\grp{\gamma_2}\cap\grp{\gamma_3,\gamma_4}$. Assume $g_2=ag_3+bg_4$. Then
    $$1=\tau(g_2,h_3)=\tau(ag_3,h_3)=\zeta^{a\alpha_{33}},$$
i.e. $a\equiv 0\pmod{\lambda}$. In the same way,
    $$1=\tau(g_2,h_4)=\zeta^{b\alpha_{44}},$$
i.e. $b\equiv 0\pmod{\lambda}$. Hence $\grp{\gamma_2}\cap\grp{\gamma_3,\gamma_4}=\{0\}$. Similarly
$\grp{\gamma_1}\cap\grp{\gamma_2,\gamma_3,\gamma_4}=\{0\}$. Hence $\grp{\gamma_i}_i=\bigoplus_i\grp{\gamma_i}$.
\end{proof}

{\samepage
From theorem~\ref{teo:p-1} we get the following probabilistic algorithm to determine gene\-rators of the
$m$-torsion subgroup $\gruppen[m]<\gruppen$, where $m\mid |\gruppen|$ is the largest divisor of $|\gruppen|$ such that
$\ell\mid p-1$ for every prime number $\ell\mid m$.

\begin{alg}\label{alg:1a}
As input we are given a hyperelliptic curve $C$ of genus two defined over a prime field $\FF_p$, the number
$N=|\gruppen|$ of $\FF_p$-rational elements of the Jacobian, and the prime factors $p_1,\dots,p_n$ of $\gcd(N,p-1)$.
The algorithm outputs elements $\gamma_i\in\gruppen[m]$ such that $\gruppen[m]=\bigoplus_i\grp{\gamma_i}$ in the
following steps.
\begin{enumerate}
 \item Set $\gamma_i:=0$, $1\leq i\leq 4$. For $\ell\in\{p_1,\dots,p_n\}$ do the following:
 \begin{enumerate}
  \item Use algorithm~\ref{alg:1} to determine elements $\tilde\gamma_i\in\gruppen_\ell$, $1\leq i\leq 4$, such that
  $\grp{\tilde\gamma_i}_i=\bigoplus_i\grp{\tilde\gamma_i}$.\label{step:kandidater}
  \item If $\prod_i|\tilde\gamma_i|<|\gruppen_\ell|$, then go to step~\ref{step:kandidater}.
  \item Set $\gamma_i:=\gamma_i+\tilde\gamma_i$, $1\leq i\leq 4$.
 \end{enumerate}
 \item Output $\gamma_1$, $\gamma_2$, $\gamma_3$ and $\gamma_4$.
\end{enumerate}
\end{alg}
}

\begin{rem}
By remark~\ref{rem:runningtime}, algorithm~\ref{alg:1a} has expected running time $O(\log p)$. Hence
algorithm~\ref{alg:1a} is an efficient, probabilistic algorithm to determine generators of the $m$-torsion subgroup
$\gruppen[m]<\gruppen$, where $m\mid |\gruppen|$ is the largest divisor of $|\gruppen|$ such that $\ell\mid p-1$ for
every prime number $\ell\mid m$.
\end{rem}

\begin{rem}
The strategy of algorithm~\ref{alg:1} can be applied to \emph{any} finite, abelian group $\gruppen$ with bilinear,
non-degenerate pairings into cyclic groups. For the strategy to be efficient, the pairings must be
efficiently computable, and the discrete logarithm problem in the cyclic groups must be easy.
\end{rem}

\end{document}